# A note on the Sylvester equation over max-plus algebra

Pingke Li

Department of Industrial Engineering, Tsinghua University, Beijing, 100084, P.R. China

**ABSTRACT**

It is known that the solvability of a Sylvester equation over max-plus algebra can be determined in polynomial time by verifying its principal solution. A succinct representation of the principal solution is presented, with a more accurate computational complexity, for a Sylvester equation over max-plus algebra.

**KEYWORDS**

max-plus algebra; Sylvester equation; Kronecker product

## 1. Introduction

Max-plus algebra is the linear algebra over the set $\mathbb{R}_{\max} = \mathbb{R} \cup \{-\infty\}$ endowed with maximization as addition and the usual sum as multiplication, i.e., $a \oplus b = \max\{a, b\}$ and $a \otimes b = a + b$ for $a, b \in \mathbb{R}_{\max}$. Analogously, min-plus algebra can be defined over $\mathbb{R}_{\min} = \mathbb{R} \cup \{+\infty\}$ with $a \oplus' b = \min\{a, b\}$ and $a \otimes' b = a + b$ for $a, b \in \mathbb{R}_{\min}$. The pair of operations $(\oplus, \otimes)$, as well as the pair $(\oplus', \otimes')$, can be extended to matrices and vectors in the same way as in linear algebra.

Max-plus algebra provides an attractive approach to describing some discrete event systems and optimization problems in production and transportation due to the possibility of handling some types of nonlinear problems in a linear-like manner. The reader may refer to Baccelli et al. (1992), Heidergott et al. (2005), and Butkovič (2010) for more information on the theory and application of max-plus algebra.

We consider in this note the following max-plus Sylvester equation

$$\bigoplus_{k=1}^{p} A_k \otimes X \otimes B_k = C \qquad (1)$$

where $A_k \in \mathbb{R}_{\max}^{m \times m}$, $B_k \in \mathbb{R}_{\max}^{n \times n}$, $k = 1, 2, \ldots, p$, and $C \in \mathbb{R}_{\max}^{m \times n}$ are given matrices while $X \in \mathbb{R}_{\max}^{m \times n}$ is an unknown matrix. Butkovič and Fiedler (2011) first investigated such type of matrix equations and demonstrated that a max-plus Sylvester equation can be transformed, by means of the max-plus Kronecker product operation, into a system of max-plus linear equations which have been intensively studied in the literature.

Hashemi, Khalilabadi, and Tavakolipour (2015) focused on a special case of the max-plus Sylvester equation $A \otimes X \oplus X \otimes B = C$ and indicated that its solvability can be verified with a time complexity of $O(m^2n + mn^2)$. Besides, Myšková (2016) addressed an interval version of the max-plus matrix equation $A \otimes X \otimes B = C$.

In this note, we demonstrate that the solvability of the max-plus Sylvester equation (1) can be verified in general with a time complexity of $O(p(m^2n + mn^2))$ taking advantage of the properties of max-plus and min-plus Kronecker products. Note that the operations of max-plus and min-plus Kronecker products are defined analogously as their counterparts in linear algebra with the multiplication replaced by the usual sum. Some properties of the max-plus and min-plus Kronecker products have been discussed by Butkovič and Fiedler (2011), Hashemi, Khalilabadi, and Tavakolipour (2015), and Hashemi, Tavakolipour, and Shirazi (2014).

The rest of this note is organized as follows. Some basic results on solving max-plus linear equations are recalled in Section 2. The solvability of a max-plus Sylvester equation is discussed in Section 3 by formulating its principal solution in a succinct manner. A short conclusion is included in Section 4.

## 2. Max-plus linear systems

One of the basic problems in max-plus algebra is that of solving a system of max-plus linear equations of the form

$$A \otimes \mathbf{x} = \mathbf{b} \qquad (2)$$

where $A \in \mathbb{R}_{\max}^{m \times n}$ and $\mathbf{b} \in \mathbb{R}_{\max}^{m}$ are given while $\mathbf{x} \in \mathbb{R}_{\max}^{n}$ is an unknown vector. Without loss of generality, we may assume that $A$ is doubly $\mathbb{R}$-astic, i.e., $A$ has in each row and each column at least one finite element and also $\mathbf{b}$ consists of only finite elements.

It is now well known that $A \otimes \mathbf{x} \leq \mathbf{b}$ if and only if $\mathbf{x} \leq A^\sharp \otimes' \mathbf{b}$ where $A^\sharp = -A^T$ is the conjugate of $A$ in the context of max-plus algebra. Consequently, it suffices to check whether or not $A^\sharp \otimes' \mathbf{b}$, called the principal solution, is indeed a solution in order to determine the solvability of $A \otimes \mathbf{x} = \mathbf{b}$.

**Lemma 2.1.** *[Butkovič (2010)] A max-plus linear system $A \otimes \mathbf{x} = \mathbf{b}$ is solvable if and only if its principal solution*

$$\mathbf{x}^* = A^\sharp \otimes' \mathbf{b} \qquad (3)$$

*is a solution to $A \otimes \mathbf{x} = \mathbf{b}$. When $\mathbf{x}^*$ is a solution, it is also the maximum solution to $A \otimes \mathbf{x} = \mathbf{b}$ with respect to the canonical order of vectors over max-plus algebra.*

Lemma 2.1 and its equivalent variants have been known for a long time since the early development of the theory of max-plus algebra in 1960s. The reader may refer to Butkovič (2010) for a detailed discussion on the role of the principal solution in solving max-plus linear systems. By Lemma 2.1, the principal solution $\mathbf{x}^* = A^\sharp \otimes' \mathbf{b}$ can be calculated and verified with a time complexity of $O(mn)$ for a given max-plus linear system $A \otimes \mathbf{x} = \mathbf{b}$.

Now we consider the problem of solving the max-plus linear matrix equation

$$A \otimes X \otimes B = C \tag{4}$$

where $A \in \mathbb{R}_{\max}^{m \times m}$, $B \in \mathbb{R}_{\max}^{n \times n}$, and $C \in \mathbb{R}_{\max}^{m \times n}$ are given matrices while $X \in \mathbb{R}_{\max}^{m \times n}$ is an unknown matrix. It can be viewed as a degenerated scenario of the max-plus Sylvester equation (1). A brute force way to solve $A \otimes X \otimes B = C$, as illustrated in Butkovič and Fiedler (2011), is to reformulate it into a max-plus linear system as

$$(B^T \boxtimes A) \otimes \text{vec}(X) = \text{vec}(C) \tag{5}$$

where the $\text{vec}(\cdot)$ operator transforms a matrix into a vector by stacking its columns one underneath the other and $\boxtimes$ denotes the max-plus Kronecker product operation, i.e.,

$$B^T \boxtimes A = \begin{bmatrix} b_{11} \otimes A & \cdots & b_{n1} \otimes A \\ \vdots & \ddots & \vdots \\ b_{1n} \otimes A & \cdots & b_{nn} \otimes A \end{bmatrix}. \tag{6}$$

Such a reformulation, analogous to its counterpart in linear algebra, holds because it requires only that the commutative, associative, and distributive laws hold for the pair of operations $(\oplus, \otimes)$. Note that the left Kronecker product form was used in Butkovič and Fiedler (2011) while in this note the right Kronecker product form is adopted.

**Lemma 2.2.** *[Butkovič and Fiedler (2011)] A max-plus linear matrix equation $A \otimes X \otimes B = C$ is solvable if and only if the principal solution $X^*$ is a solution where*

$$\text{vec}(X^*) = (B^T \boxtimes A)^\sharp \otimes' \text{vec}(C). \tag{7}$$

It is clear that Lemma 2.2 is a direct consequence of Lemma 2.1. Accordingly, the solvability of $A \otimes X \otimes B = C$ can be determined with a polynomial time complexity of $O(m^2 n^2)$ if a direct calculation procedure is performed. However, taking advantage of the properties of max-plus and min-plus Kronecker products, it turns out that the principal solution $X^*$ to $A \otimes X \otimes B = C$ can be calculated in a more efficient manner.

Note that $(B^T \boxtimes A)^\sharp = (B^T)^\sharp \boxtimes' A^\sharp = (B^\sharp)^T \boxtimes' A^\sharp$ where $\boxtimes'$ denotes the min-plus

Kronecker product operation. Consequently,

$$\text{vec}(X^*) = (B^T \boxtimes A)^\sharp \otimes' \text{vec}(C) = ((B^\sharp)^T \boxtimes' A^\sharp) \otimes' \text{vec}(C) \tag{8}$$

which is actually the reformulation of

$$X^* = A^\sharp \otimes' C \otimes' B^\sharp \tag{9}$$

via the min-plus Kronecker product operation. Therefore, it takes only a time complexity of $O(m^2n + mn^2)$ to calculate the principal solution $X^*$ and check the solvability of $A \otimes X \otimes B = C$.

Although this result can be derived in a direct manner according to the properties of residuated functions as discussed in Baccelli et al. (1992), the use of the Kronecker product operations offers an alternative and succinct approach which can be extended readily to the general max-plus Sylvester equations. This result can also be applied to deal with the interval version of $A \otimes X \otimes B = C$ considered in Myšková (2016).

## 3. Max-plus Sylvester equations

The basic type of max-plus Sylvester equations is of the form

$$A_1 \otimes X \otimes B_1 \oplus A_2 \otimes X \otimes B_2 = C \tag{10}$$

where $A_1, A_2 \in \mathbb{R}_{\max}^{m \times m}$, $B_1, B_2 \in \mathbb{R}_{\max}^{n \times n}$, and $C \in \mathbb{R}_{\max}^{m \times n}$ are given matrices while $X \in \mathbb{R}_{\max}^{m \times n}$ is an unknown matrix. It can be reformulated as

$$(B_1^T \boxtimes A_1 \oplus B_2^T \boxtimes A_2) \otimes \text{vec}(X) = \text{vec}(C) \tag{11}$$

and hence, its principal solution $X^*$ can be calculated as

$$\begin{aligned}
\text{vec}(X^*) &= (B_1^T \boxtimes A_1 \oplus B_2^T \boxtimes A_2)^\sharp \otimes' \text{vec}(C) \\
&= ((B_1^T \boxtimes A_1)^\sharp \oplus' (B_2^T \boxtimes A_2)^\sharp) \otimes' \text{vec}(C) \\
&= ((B_1^\sharp)^T \boxtimes A_1^\sharp \oplus' (B_2^\sharp)^T \boxtimes A_2^\sharp) \otimes' \text{vec}(C) \\
&= ((B_1^\sharp)^T \boxtimes A_1^\sharp) \otimes' \text{vec}(C) \oplus' ((B_2^\sharp)^T \boxtimes A_2^\sharp) \otimes' \text{vec}(C). \tag{12}
\end{aligned}$$

This fact leads to the following solvability criterion for a max-plus Sylvester equation of this basic type.

**Theorem 3.1.** *A max-plus Sylvester equation $A_1 \otimes X \otimes B_1 \oplus A_2 \otimes X \otimes B_2 = C$ is*

*solvable if and only if the principal solution $X^*$ is a solution where*

$$X^* = (A_1^\sharp \otimes' C \otimes' B_1^\sharp) \oplus' (A_2^\sharp \otimes' C \otimes' B_2^\sharp). \tag{13}$$

As a consequence of Theorem 3.1, the solvability of $A_1 \otimes X \otimes B_1 \oplus A_2 \otimes X \otimes B_2 = C$ can be verified with a time complicity of $O(m^2n + mn^2)$. Note that the equation considered by Hashemi, Khalilabadi, and Tavakolipour (2015) is a special scenario of this type by setting $B_1$ and $A_2$ the max-plus unit matrices of compatible dimensions, i.e., the square matrices whose diagonal elements are all zero and off-diagonal elements are all $-\infty$.

Moreover, it is clear that Theorem 3.1 can be extended to general max-plus Sylvester equations.

**Theorem 3.2.** *A max-plus Sylvester equation*

$$\bigoplus_{k=1}^{p} A_k \otimes X \otimes B_k = C \tag{14}$$

*is solvable if and only if the principal solution $X^*$ is a solution where*

$$X^* = \bigoplus_{k=1}^{p}{}' A_k^\sharp \otimes' C \otimes' B_k^\sharp \tag{15}$$

*which can be calculated and verified with a time complexity of $O(p(m^2n + mn^2))$.*

Note that by transforming it into a max-plus linear system, the solvability criterion of a max-plus Sylvester equation was first presented in Butkovič and Fiedler (2011) using the left Kronecker product form, which may result in Theorem 3.2 applying the properties of max-plus and min-plus Kronecker products. Besides, Theorem 3.2 can be derived directly as well within the framework of residuation theory.

## 4. Conclusions

Following the results in Butkovič and Fiedler (2011), the principal solution of a Sylvester equation over max-plus algebra is presented in a succinct manner in this note. It reveals a more accurate computational complexity for determining the solvability of a max-plus Sylvester equation and also generalizes the results in Hashemi, Khalilabadi, and Tavakolipour (2015). The presented approach can be extended as well to Sylvester equations defined over some analogous semiring structures such as fuzzy algebra.

**Funding**

This work was supported by the National Natural Science Foundation of China under Grant 61203131 and the Tsinghua University Initiative Scientific Research Program under Grant 2014z21017.